\documentclass[12pt]{article}
\usepackage{amsmath,amssymb}
\usepackage{graphicx}
\newcommand{\eqdef}{\stackrel{\text{def}}{=}}

\newcommand{\n}{\nonumber \\}
\newcommand{\bm}{\boldsymbol}
\newcommand{\ignore}[1]{}
\numberwithin{equation}{section}
\newcommand{\Romannumeral}[1]{\uppercase\expandafter{\romannumeral#1}}

\newtheorem{theo}{\bf Theorem}[section]

\newtheorem{rema}[theo]{\bf Remark}

\newtheorem{prop}[theo]{\bf Proposition}

\allowdisplaybreaks[4]

\begin{document}

\baselineskip=20pt
\newcommand{\preprint}{
\vspace*{-20mm}\begin{flushleft}\end{flushleft}
}
\newcommand{\Title}[1]{{\baselineskip=26pt
  \begin{center} \Large \bf #1 \\ \ \\ \end{center}}}
\newcommand{\Author}{\begin{center}
  \large \bf 
  Ryu Sasaki${}$ \end{center}}
\newcommand{\Address}{\begin{center}
     Department of Physics and Astronomy, Tokyo University of Science,
     Noda 278-8510, Japan
        \end{center}}
\newcommand{\Accepted}[1]{\begin{center}
  {\large \sf #1}\\ \vspace{1mm}{\small \sf Accepted for Publication}
  \end{center}}

\preprint
\thispagestyle{empty}

\Title{Rahman polynomials}

\Author

\Address
\vspace{1cm}

\begin{abstract}
Two very closely related  Rahman polynomials are constructed explicitly as the left eigenvectors of 
certain multi-dimensional discrete time Markov chain operators 
$K_n^{(i)}({\boldsymbol x},{\boldsymbol y};N)$, $i=1,2$.
They are convolutions of an $n+1$-nomial distribution $W_n({\boldsymbol x};N)$ and an $n$-tuple of 
binomial distributions $\prod_{i}W_1(x_i;N)$. The one for the original Rahman polynomials is
$K_n^{(1)}({\boldsymbol x},{\boldsymbol y};N)
=\sum_{\boldsymbol z}W_n({\boldsymbol x}-{\boldsymbol z};N-\sum_{i}z_i)
\prod_{i}W_1(z_i;y_i)$. The closely related one is
\  $K_n^{(2)}({\boldsymbol x},{\boldsymbol y};N)
=\sum_{\boldsymbol z}W_n({\boldsymbol x}-{\boldsymbol z};N-\sum_{i}y_i)
\prod_{i}W_1(z_i;y_i)$.
The original Markov chain was introduced and discussed by Hoare, Rahman and Gr\"{u}nbaum
as a multivariable version of the known soluble single variable one.
The new one is a generalisation of that of Odake and myself.
The anticipated solubility of the model gave Rahman polynomials the prospect of the first
multivariate hypergeometric function of Aomoto-Gelfand type connected with solvable dynamics. 
The promise is now realised. The  $n^2$ system parameters $\{u_{i\,j}\}$ of
the Rahman polynomials are completely determined.
These $u_{i\,j}$'s are irrational functions  of the original system parameters, 
the probabilities of the multinomial and binomial distributions.
\end{abstract}

%
%
\section{Introduction}
\label{sec:intro}

The two types of Rahman polynomials,  typical examples of multivariate orthogonal polynomials,
\cite{coo-hoa-rah77, HR0, HR, gr1, gr2, gr3, IT, I}, are constructed explicitly.
The method of construction is very different from the existing theories and models of multivariate
orthogonal polynomials, see for example, 
\cite{Gri1, Gri10, dunkl0, tra, zheda, dunkl, mizu1, mt, HR, gr1, ilxu,  gr2, IT, mizu, gr3, I, I11, miki1,
diaconis13, xu, genest0, genest1, genest2, genest3, shibukawa, Gri11, Gri2, genest4, miki2}.
The Rahman polynomials and the two multivariate orthogonal polynomials constructed recently, 
the Krawtchouk 
\cite{I23, mKrawt} and Meixner polynomials \cite{mMeix} are a certain type of 
Aomoto-Gelfand \cite{AK, gelfand} terminating hypergeometric functions, see \eqref{Pm}.
Like other orthogonal polynomials of a discrete variable \cite{askey, ismail, koeswart, os12}, 
the knowledge of $n$ degree one polynomials is enough to determine the system parameters 
$\{u_{i\,j}\}$ \eqref{u1def}, \eqref{u2def} completely. 
The entire polynomials are expressed by $u_{i\,j}$'s polynomially, see  formula \eqref{Pm}.
A certain dynamical setting, a quite essential one, determines the degree one polynomials.
For the Rahman polynomials, it is the discrete time Markov chain and for the multivariate
Krawtchouk and Meixner polynomials, it is the birth and death processes
 \cite{feller, bdsol, mKrawt,mMeix}.
These dynamical settings determine the orthogonal weight functions at the same time.
But the polynomials are determined independently of the orthogonality weights.

The construction of the Rahman polynomials is straightforward.
The problem setting was done completely by the works of Hoare, Rahman and Gr\"{u}nbaum 
\cite{HR, gr1, gr2, gr3}. 
I only have to calculate the degree one {\em left eigenpolynomials} of the two Markov chain operators
$K_n^{(i)}({\bm x},{\bm y}; N,{\bm\alpha}, {\bm \beta})$, $i=1,2$ and determine the system parameters 
$\{u_{i\,j}\}$ as  functions of $2n$ system parameters ${\bm \alpha}$ and ${\bm \beta}$.
The remaining jobs are to determine the entire spectrum 
$\mathcal{E}({\bm m})=\prod_{i=1}^n\lambda_i^{m_i}$ \eqref{emformula}, in which $\lambda_i$'s
are the eigenvalues of the degree one polynomials.
The final task is to demonstrate that the entire polynomials are the {\em left eigenvectors} of the
Markov chain operators with the multiplicative spectrum. The generating function by Mizukawa \cite{mizu}
plays the central role.

 This paper is organised as follows. 
 The problem setting of discrete time Markov chains due to Hoare, Rahman and Gr\"{u}nbaum is
 briefly recapitulated in section \ref{sec:set}.
 The important role played by left eigenvectors are stated in {\bf Theorem \ref{lefttheo}}.
 The main players, Aomoto-Gelfand \cite{AK,gelfand} terminating hypergeometric functions and the generating function
 are introduced in {\bf Theorem \ref{theo:pm}}.
 The two Rahman polynomials are constructed in section \ref{sec:nRah}.
 Starting with  degree one left eigenpolynomials in \S\ref{sec:degone}, 
 the spectrum is determined in \S\ref{sec:eigval} and the proof of the left eigenvector equations 
 for $P_{\bm m}({\bm x};{\bm u})$ is provided in \S\ref{sec:lefteigeq}.
 The final section is for some comments.
This paper is prepared plainly so that non-experts can understand.

\section{Problem setting}
\label{sec:set}

Let us start with a brief review of the general framework of  
discrete time Markov chains of $n$-variables.
For simplicity of presentation, positive transition matrices on finite sets only are discussed. 
\subsection{Discrete time Markov chains of $n$-variables}
\label{sec:genMarkov}
The main actor is the transition matrix $\mathcal{K}_n({\bm x},{\bm y})$ which governs the 
evolution of the probability distribution function $\mathcal{P}({\bm x};\ell)$ over $\Xi$, a finite set,
\begin{equation*}
\mathcal{P}({\bm x};\ell)\longrightarrow \mathcal{P}({\bm x};\ell+1)
=\sum_{{\bm y}\in\Xi}\mathcal{K}_n({\bm x},{\bm y})\mathcal{P}({\bm y};\ell),
\quad \mathcal{P}({\bm x};\ell)\ge0,
\quad \sum_{{\bm x}\in\Xi}\mathcal{P}({\bm x};\ell)=1,\  \ell\in\mathbb{N}_0, 
\end{equation*}
at time step $\ell$ to $\ell+1$.
That is, $\mathcal{K}_n({\bm x},{\bm y})$ specifies
 the one step transition probability from point
${\bm y}$ to ${\bm x}$ which belong to a finite set $\Xi$,  a subset  of the non-negative integer lattice 
in $n$-dimensions $\mathbb{N}_0^n$,
\begin{align*} 
\mathcal{K}_n({\bm x},{\bm y})>0,\quad {\bm x}=(x_1,\ldots,x_n),
\quad {\bm y}=(y_1,\ldots,y_n),\quad 
\qquad {\bm x},{\bm y}\in\Xi\subseteq\mathbb{N}_0^n,
\end{align*}
satisfying
\begin{align}
&\text{\bf conservation of probability:}
\quad  \sum_{{\bm x}\in\Xi}\mathcal{K}_n({\bm x},{\bm y})=1,
\label{conspr}\\
&\qquad \qquad \Longrightarrow 
\sum_{{\bm x}\in\Xi}\mathcal{P}({\bm x};\ell+1)=\sum_{{\bm x}\in\Xi}\mathcal{P}({\bm x};\ell)=1.
\nonumber
\end{align}
Let us impose the condition that $\mathcal{P}({\bm x};\ell)$ does not explode and 
approaches to a {\em stationary} distribution $\pi_S$ at large $\ell$;
\begin{equation*}
\lim_{\ell\to\infty}\mathcal{P}({\bm x};\ell)=\pi_S({\bm x}).
\end{equation*}
The existence of the stationary distribution is guaranteed  
when $\mathcal{K}_n({\bm x},{\bm y})$
satisfies the following {\em reversibility} or {\em detailed balance} condition
with  a certain positive {\em reversible distribution} $\pi_R({\bm x})$,
\begin{align}
\text{\bf reversible:}\qquad & \mathcal{K}_n({\bm x},{\bm y})\pi_R({\bm y})
=\mathcal{K}_n({\bm y},{\bm x})\pi_R({\bm x}),\quad \pi_R({\bm x})>0,
\quad {\bm x}, {\bm y}\in\Xi.
\label{revcond}
\end{align}
Deviding the above reversibility condition by $\sqrt{\pi_R({\bm x})\pi_R({\bm y})}$, 
one obtains a real symmetric matrix $\mathcal{T}_n$;
\begin{align}
\mathcal{T}_n({\bm x},{\bm y})&\eqdef
\frac1{\sqrt{\pi_R({\bm x})}}\,\mathcal{K}_n({\bm x},{\bm y})\sqrt{\pi_R({\bm y})}
=\frac1{\sqrt{\pi_R({\bm y})}}\,\mathcal{K}_n({\bm y},{\bm x})\sqrt{\pi_R({\bm x})}
\label{Tdef}\\[2pt]
&=\mathcal{T}_n({\bm y},{\bm x}),\qquad {\bm x}, {\bm y}\in\Xi.
\label{symm}
\end{align}
In other words, the  transition matrix $\mathcal{K}_n({\bm x},{\bm y})$ is related to a
real positive symmetric matrix $\mathcal{T}_n({\bm x},{\bm y})$ by a similarity transformation 
in terms of $\sqrt{\pi_R({\bm x})}$.
This means the following
\begin{prop}
\label{eig1}
The eigenvalues of $\mathcal{K}_n$ are real and the range of the spectrum is
\begin{equation}
  -1<\text{The  eigenvalues of }\mathcal{K}_n({\bm x},{\bm y})\le1,
\end{equation}
and the  eigenvector of the maximal eigenvalue 1 is the 
the reversible  distribution 
\begin{align}
  \sum_{{\bm y}\in\Xi}\mathcal{K}_n({\bm x},{\bm y})\pi_R({\bm y})
=\pi_R({\bm x})\sum_{{\bm y}\in\Xi}\mathcal{K}_n({\bm y},{\bm x})=\pi_R({\bm x}).
\label{maxeig}
\end{align}
The maximal eigenvalue is simple.
These are due to Perron's theorem on positive matrices.
As $\ell$ becomes large, all the components of $\mathcal{P}({\bm x};\ell)$ 
belonging to the eigenvalues other than 1 diminish and
\begin{equation*}
\lim_{\ell\to\infty}\mathcal{P}({\bm x};\ell)=\pi_R({\bm x}).
\end{equation*}
In other words, the reversible distribution is stationary.
\end{prop}
The above maximal eigenvector relation \eqref{maxeig} also means the following
\begin{rema}
\label{unitleft}
A constant is the left eigenvector of $\mathcal{K}_n({\bm x},{\bm y})$ with the maximal
eigenvalue 1;
\begin{equation*}
v_{\bm 0}({\bm x})\equiv1,\quad
\sum_{{\bm x}\in\Xi}\mathcal{K}_n({\bm x},{\bm y})v_{\bm 0}({\bm x})=
\sum_{{\bm x}\in\Xi}\mathcal{K}_n({\bm x},{\bm y})=1=v_{\bm 0}({\bm y}).
\end{equation*}
\end{rema}
By generalising the {\bf Remark \ref{unitleft}}, one arrives at the following
\begin{theo}
\label{lefttheo}
The complete set of left eigenvectors $\{v_{\bm m}({\bm x})\}$ of 
$\mathcal{K}_n({\bm x},{\bm y})$, 
\begin{equation}
\sum_{{\bm x}\in\Xi}\mathcal{K}_n({\bm x},{\bm y})v_{\bm m}({\bm x})
=\kappa({\bm m})v_{\bm m}({\bm y}),\quad -1<\kappa({\bm m})\leq1=\kappa({\bm 0}),\quad 
\forall{\bm m}\in\Xi^d,
\label{leigeq}
\end{equation}
provides the complete sets of eigenvectors of $\mathcal{K}_n({\bm x},{\bm y})$
and $\mathcal{T}_n({\bm x},{\bm y})$;
\begin{align}
&\sum_{{\bm y}\in\Xi}\mathcal{K}_n({\bm x},{\bm y})v_{\bm m}({\bm y})\pi_R({\bm y})
=\kappa({\bm m})v_{\bm m}({\bm x})\pi_R({\bm x}),\quad \forall{\bm m}\in\Xi^d,
\label{Kneig1}\\
&\sum_{{\bm y}\in\Xi}\mathcal{T}_n({\bm x},{\bm y})v_{\bm m}({\bm y})\sqrt{\pi_R({\bm y})}
=\kappa({\bm m})v_{\bm m}({\bm x})\sqrt{\pi_R({\bm x})},\quad \forall{\bm m}\in\Xi^d.
\label{Tneig1}
\end{align}
Here  $\Xi^d\subseteq\mathbb{N}_0^n$ is the dual set of $\Xi$ and $\#(\Xi^d)=\#(\Xi)$.
The left eigenvectors are orthogonal to each other with the orthogonality weight
function $\pi_R({\bm x})$;
\begin{equation}
\sum_{{\bm x}\in\Xi}v_{\bm m}({\bm x})v_{{\bm m}'}({\bm x})\pi_R({\bm x})=0,\quad 
{\bm m}\neq{\bm m}',
\label{vmort}
\end{equation}
since $\mathcal{K}_n({\bm x},{\bm y})$ has some parameters and for generic values of the parameters, 
the eigenvalues are not degenerate.
\end{theo}

This setup of Markov chains in one dimensions has been discussed by many authors 
\cite{coo-hoa-rah77, HR0, os39}, in connection with ``cumulative Bernoulli trials."
The multivariate version was initiated by Hoare, Rahman and Gr\"unbaum \cite{HR, gr1,gr2,gr3}
in a quest for a new multivariate version of the Krawtchouk polynomials, called Rahman polynomials.
In particular, the explicit setup for $K_n^{(1)}$  in the next subsection was due to 
Gr\"unbaum  and Rahman \cite{gr3}, and this paper is cite as I.
\subsection{Quest for Rahman polynomials}
\label{sec:quest}

The purpose of this paper is to present the explicit forms of two types of  Rahman polynomials 
by solving two  $n$-variable Markov chains with the
transition matrix (operator) $K_n^{(i)}$, $i=1,2$ composed of the
convolution of an $n+1$-nomial distribution $W_n$  and $n$ binomial distributions $W_1$;
\begin{align} 
  & K_n^{(1)}({\bm x},{\bm y}; N,{\bm\alpha}, {\bm \beta})
 \eqdef\sum_{{\bm z}\in\mathcal{X}}W_n({\bm x}-{\bm z};N-|z|,{\bm \beta})
 \prod_{i=1}^nW_1(z_i;y_i,\alpha_i)>0,
\label{Kn1def}\\
 & K_n^{(2)}({\bm x},{\bm y}; N,{\bm\alpha}, {\bm \beta})
 \eqdef\sum_{{\bm z}\in\mathcal{X}}W_n({\bm x}-{\bm z};N-|y|,{\bm \beta})
 \prod_{i=1}^nW_1(z_i;y_i,\alpha_i)>0,
\label{Kn2def}
\end{align}
in which
\begin{align*} 
&{\bm x}=(x_1,\ldots,x_n),\quad {\bm y}=(y_1,\ldots,y_n),\quad {\bm z}=(z_1,\ldots,z_n), \quad 
|x|\eqdef\sum_{i=1}^nx_i, \quad N\in\mathbb{N},
\\
&\qquad {\bm x},{\bm y},{\bm z}\in\mathcal{X}\eqdef\{\bm{x}\in\mathbb{N}_0^n\ |\, |x|\le N\},
\qquad \quad \#\bigl(\mathcal{X}\bigr)=\binom{N+n}{n},\n
&{\bm \alpha}=(\alpha_1,\ldots,\alpha_n),\ {\bm \beta}=(\beta_1,\ldots,\beta_n),
\  0<\alpha_i,\beta_i<1,\quad i=1,\ldots,n,\quad 0<|\beta|<1,
\end{align*}
and  the  binomial  $W_1$ and multinomial $W_n$ distributions are,
\begin{align}
&W_1(x;y,\alpha)\eqdef \binom{y}{x}\alpha^x(1-\alpha)^{y-x}>0,
\quad \sum_{x=0}^yW_1(x;y,\alpha)=1,
\label{W1def}\\
&W_n({\bm x};N,{\bm \beta})\eqdef 
\frac{N!\cdot(1-|\beta|)^{N-|x|}}{x_1!\cdots x_n!(N-|x|)!}\cdot\prod_{i=1}^n\beta_i^{x_i}
=\binom{N}{\bm{x}}\beta_0^{x_0}\bm{\beta}^{\bm{x}}>0,
\label{Wndef}\\
&\sum_{{\bm x}\in\mathcal{X}}W_n({\bm x};N,{\bm \beta})=1,\quad    x_0\eqdef
N-|x|,\  \beta_0\eqdef 1-|\beta|,\  \binom{N}{\bm{x}}\eqdef\frac{N!}{x_1!\cdots x_n!x_0!}.
\nonumber
\end{align}
To be more explicit, the $z_i$ summation in the  Markov chains \eqref{Kn1def}, \eqref{Kn2def} are
\begin{align*}
&\text{for} \quad K_n^{(1)} \quad 0\leq z_i\leq{\text{min}}(x_i,y_i), \qquad i=1,\ldots,n,\\
&\text{for} \quad K_n^{(2)}  \quad 
\text{max}(0,x_i+|y|-N)\leq z_i\leq \text{min}(x_i,y_i),\qquad i=1,\ldots,n,\\
&\Longrightarrow \text{conservation of probability}\quad 
\sum_{{\bm x}\in\mathcal{X}} K_n^{(i)}({\bm x},{\bm y}; N,{\bm\alpha}, {\bm \beta})=1,\quad i=1,2.
\end{align*}
The conservation of probability 
holds as the double sum $\sum_{{\bm x},{\bm z}\in\mathcal{X}}$ can be changed.
Summing over $\bm x$ for fixed $\bm z$ gives 
$\sum_{{\bm x}\in{\mathcal X}}W_n({\bm x}-{\bm z};N-|z|,{\bm \beta})=1$,
$\sum_{{\bm x}\in{\mathcal X}}W_n({\bm x}-{\bm z};N-|y|,{\bm \beta})=1$.
The remaining $\bm z$ sum over $W_1$'s gives trivially 1.

\bigskip
The reversible distribution for $K_n^{(1)}$ was identified in I \cite{gr3}. 
It is the same multinomial distribution as 
$W_n({\bm x};N,{\bm \beta})$ but with a different set of probability parameters ${\bm \eta}$,
\begin{align}
& K_n^{(1)}({\bm x},{\bm y}; N,{\bm\alpha}, {\bm \beta})W_n({\bm y};N,{\bm \eta})
 = K_n^{(1)}({\bm y},{\bm x}; N,{\bm\alpha}, {\bm \beta})W_n({\bm x};N,{\bm \eta}),
 \label{Kn1rev}\\
& \eta_i=\frac{\beta_i}{1-\alpha_i}\frac1{D_n},\quad i=1,\ldots,n,
\qquad D_n\eqdef 1+\sum_{k=1}^n\frac{\alpha_k\beta_k}{1-\alpha_k},
\label{etadef1}\\
&\frac{1-\alpha_1}{\beta_1}\eta_1=\frac{1-\alpha_2}{\beta_2}\eta_2=\cdots=
\frac{1-\alpha_n}{\beta_n}\eta_n=\frac{1-|\eta|}{1-|\beta|}=
1-\sum_{k=1}^n\alpha_k\eta_k=D_n^{-1}.
\tag{I.1.6}
\end{align}
For $n=1$, $\eta=\beta/(1-\alpha+\alpha\beta)$ as reported in I \cite{gr3}.
I will not repeat the derivation of these complicated formulas (I.1.6).
Instead I show a path to \eqref{etadef1} and (I.1.6) by using an equivalent expression 
of the multinomial distribution \eqref{Wndef}, which derives from the multinomial
expansion
\begin{equation*}
\left(1+\sum_{i=1}^n q_i\right)^N
=\sum_{{\bm x}\in\mathcal{X}}\binom{N}{{\bm x}}\prod_{i=1}^nq_i^{xi},\quad q_i\in\mathbb{C},
\quad i=1,\ldots,n.
\end{equation*}
This leads to another form of the multinomial distribution
\begin{align}
\bar{W}_n({\bm x};N,{\bm q})\eqdef\binom{N}{{\bm x}}\prod_{i=1}^nq_i^{x_i}(1+|q|)^{-N}
&\Leftrightarrow
W_n({\bm x};N,{\bm \beta})=\binom{N}{{\bm x}}\prod_{i=1}^n\beta_i^{x_i}(1-|\beta|)^{N-|x|},\n
0<q_i=\frac{\beta_i}{1-|\beta|} &\Leftrightarrow \beta_i=\frac{q_i}{1+|q|},  \quad i=1,\ldots,n.
\label{qbcorres}
\end{align}
This type of multinomial distributions was useful for the derivation of the multivariate
Krawtchouk and Meixner  polynomials based on the Birth and Death formalism \cite{mKrawt,mMeix}.
In terms of these, the left hand side of the reversibility condition \eqref{Kn1rev} reads
\begin{align*}
&\sum_{{\bm z}\in\mathcal{x}}\bar{W}_n({\bm x}-{\bm z};N-|z|,{\bm q})
\prod_{i=1}^n\bar{W}_1(z_i;y_i,p_i)\bar{W}_n({\bm y};N,{\bm r})\\
&=\prod_{i=1}^n\left[q_i^{x_i}\left(\frac{r_i}{1+p_i}\right)^{y_i}\right]
\sum_{\bm z}\frac{N!(N-|z|)!\left((1+|q|)(1+|r|)\right)^{-N}}
{(N-|x|)!(N-|y|)!z_i!(x_i-z_i)!(y_i-z_i)!}.
\end{align*}
The reversibility requires that this expression is ${\bm x}\leftrightarrow{\bm y}$ symmetric.
It is necessary and sufficient that the unknown parameters $\{r_i\}$ should be
\begin{equation}
r_i=q_i(1+p_i),\quad i=1,\ldots,n,
\end{equation}
which correspond to the known parameters $\{\eta_i\}$ in \eqref{etadef1} 
by the following translation rule together with \eqref{qbcorres},
\begin{equation}
r_i=\frac{\eta_i}{1-|\eta|} \Leftrightarrow \eta_i=\frac{r_i}{1+|r|}, \quad 
p_i=\frac{\alpha_i}{1-\alpha_i}\Leftrightarrow \alpha_i=\frac{p_i}{1+p_i}
\quad i=1,\ldots,n.
\label{retacorres}
\end{equation}
This might give some idea how the detailed balance condition works.\\
The reversible distribution  for the second Markov chain is
\begin{align}
& K_n^{(2)}({\bm x},{\bm y}; N,{\bm\alpha}, {\bm \beta})W_n({\bm y};N,{\bm \eta})
 = K_n^{(2)}({\bm y},{\bm x}; N,{\bm\alpha}, {\bm \beta})W_n({\bm x};N,{\bm \eta}),
 \label{Kn2rev}\\
& \eta_i=\frac{\beta_i}{1-\alpha_i}\frac1{D_n},\quad i=1,\ldots,n,
\qquad D_n\eqdef 1+\sum_{k=1}^n\frac{\beta_k}{1-\alpha_k}.
\label{etadef2}
\end{align}
For $n=1$, $\eta=\beta/(1-\alpha+\beta)$, which was reported in \cite{os39}(4.30) for
the corresponding single variable Markov chain.
By using another form of the multinomial distribution \eqref{qbcorres}
 the left hand side of the reversibility condition \eqref{Kn2rev} reads
\begin{align*}
&\sum_{{\bm z}\in\mathcal{X}}\bar{W}_n({\bm x}-{\bm z};N-|y|,{\bm q})
\prod_{i=1}^n\bar{W}_1(z_i;y_i,p_i)\bar{W}_n({\bm y};N,{\bm r})\\
&=\prod_{i=1}^n\left[q_i^{x_i}\left(\frac{r_i(1+|q|)}{1+p_i}\right)^{y_i}\right]
\sum_{{\bm z}\in\mathcal{X}}\frac{N!\left((1+|q|)(1+|r|)\right)^{-N}(p_i/q_i)^{z_i}}
{(N-|x|-|y|+|z|)!z_i!(x_i-z_i)!(y_i-z_i)!}.
\end{align*}
It is necessary and sufficient that the unknown parameters $\{r_i\}$ should be
\begin{align*}
r_i=\frac{q_i(1+p_i)}{1+|q|}=\frac{\beta_i}{1-\alpha_i}
\Longrightarrow \eta_i&=\frac{r_i}{1+|r|}=\frac{\beta_i}{1-\alpha_i}\frac1{D_n},
\tag{\ref{etadef2}}\\
&  D_n=1+|r|= 1+\sum_{k=1}^n\frac{\beta_k}{1-\alpha_k},
\quad i=1,\ldots,n.
\end{align*}

The determination of the reversible distributions \eqref{etadef1}, \eqref{etadef2} leads to the following
\begin{theo}
\label{theo:pm}
The  multivariate polynomials orthogonal with respect to the reversible distributions 
\eqref{etadef1}, \eqref{etadef2}
are $(n+1,2n+2)$ type terminating hypergeometric function of 
Aomoto-Gelfand {\rm \cite{AK,gelfand, mizu}},
\begin{gather}
\label{Pm}
P_{\bm{m}}(\bm{x};{\bm u})
\eqdef \sum_{\substack{\sum_{i,j}c_{ij}\leq N\\
(c_{ij})\in M_{n}({\mathbb N_{0}})}}
\frac{\prod\limits_{i=1}^{n}(-x_{i})_{\sum\limits_{j=1}^{n}c_{ij}}
\prod\limits_{j=1}^{n}(-m_{j})_{\sum\limits_{i=1}^{n}c_{ij}}}
{(-N)_{\sum_{i,j}c_{ij}}} \; \frac{\prod(u_{ij})^{c_{ij}}}{\prod c_{ij}!},\\[2pt]
\sum_{\bm{x}\in\mathcal{X}}W_n({\bm x};N,{\bm \eta})
P_{\bm{m}}(\bm{x};{\bm u})P_{\bm{m}'}(\bm{x};{\bm u})=0,\quad {\bm m}\neq{\bm m}',
\label{pmort}
\end{gather}
 in which the $n^2$ unspecified parameters $u_{i\,j}$'s are to be determined by the condition that
 they are the left eigenvectors of the transition matrix 
 $K_n^{(i)}({\bm x},{\bm y}; N,{\bm\alpha}, {\bm \beta})$, $i=1,2$, \eqref{Kn1def}, \eqref{Kn2def}.
 Their explicit expressions and the so-called orthogonality conditions will be given in 
 {\bf Theorem \ref{uformula}}  \eqref{u1def},  \eqref{u2def}, \eqref{P1ort}.
 In the formula, $(a)_n$ is the shifted factorial defined for 
 $a\in\mathbb{C}$ and nonnegative integer $n$,
$(a)_0=1$, $(a)_n=\prod_{k=0}^{n-1}(a+k)$, $n\ge1$ and
$M_{n}({\mathbb N}_{0})$ is the set of square matrices of degree $n$ with nonnegative integer
elements. 
 The generating function of the above hypergeometric function 
 \eqref{Pm} is well-known  {\rm\cite{mizu}},
\begin{align}
G(\bm{x};{\bm u}, {\bm t})&\eqdef \prod_{i=0}^{n}\left(\sum_{j=0}^{n}b_{ij}t_{j}\right)^{x_{i}}
=\sum_{\bm{m} \in\mathcal{X}}
\binom{N}{\bm{m}}P_{\bm{m}}(\bm{x};{\bm u})\bm{t}^{\bm{m}},
\label{gen}\\
&\bm{t}^{\bm{m}}\eqdef\prod_{j=1}^nt_j^{m_j},\quad t_0\eqdef1.
\end{align}
The parameters $\{b_{i\,j}\}$ are related to $\{u_{i\,j}\}$ 
\begin{align}
b_{0\,j}=b_{i\,0}=1 \ {\rm for}  \ 0\leq i,j\leq n \  {\rm and} \ b_{i\,j}=1-u_{i\,j}\  {\rm for} \ i,j=1,\ldots,n.
\label{aurel}
\end{align}
\end{theo}
By requiring that $P_{\bm{m}}(\bm{x};{\bm u})$ should obey difference equations 
of Birth and Death type 
\cite{ismail,bdsol},
another type of multivariate Krawtchouk polynomials is obtained \cite{mKrawt}.

%
%
\section{Rahman polynomials}
\label{sec:nRah}
In this section I show the explicit forms of $P_{\bm{m}}(\bm{x};{\bm u})$ \eqref{Pm}, that is,
the expression of ${\bm u}=\{u_{i\,j}\}$ in {\bf Theorem \ref{uformula}}.
Based on that, the general eigenvalues $\mathcal{E}({\bm m})$ are determined in 
{\bf Proposition \ref{prop:emformula}}.
Finally, the left eigenvalue equations
\begin{equation*}
\sum_{{\bm x}\in\mathcal{X}}K_n^{(i)}({\bm x},{\bm y}; N,{\bm\alpha}, {\bm \beta})
P_{\bm{m}}(\bm{x};{\bm u})=\mathcal{E}({\bm m})P_{\bm{m}}(\bm{y};{\bm u}),
\quad i=1,2,\quad  \forall{\bm m}\in\mathcal{X},
\end{equation*}
are proven in {\bf Theorem \ref{theo:main}}.

To begin with, the structure of the hypergeometric function 
$P_{\bm{m}}(\bm{x};{\bm u})$ \eqref{Pm}
means the following
\begin{prop}
\label{deg1enough}
The knowledge of $n$ degree one polynomials $P_{|m|=1}(\bm{x};{\bm u})$ is enough to determine 
${\bm u}=\{u_{i\,j}\}$ completely.
\end{prop}

\subsection{Degree one left eigenpolynomials}
\label{sec:degone}
In order to solve the degree one left eigenvalue problem of 
$K_n^{(i)}({\bm x},{\bm y}; N,{\bm\alpha}, {\bm \beta})$ \eqref{Kn1def}, \eqref{Kn2def}
the following summation formulas are useful,
\begin{align}
\text{\rm type\ (1)}\qquad 
\sum_{{\bm x}\in\mathcal{X}}x_iK_n^{(1)}({\bm x},{\bm y}; N,{\bm\alpha}, {\bm \beta})
&=\beta_iN-\beta_i\sum_{j=1}^n\alpha_jy_j+\alpha_iy_i,\quad i=1,\ldots,n,
\label{sum1form}\\
\text{\rm type\ (2)}\qquad 
\sum_{{\bm x}\in\mathcal{X}}x_iK_n^{(2)}({\bm x},{\bm y}; N,{\bm\alpha}, {\bm \beta})
&=\beta_i(N-|y|)+\alpha_iy_i,\quad \quad i=1,\ldots,n.
\label{sum2form}
\end{align}
These are obtained by combining the recurrence relations of $W_1$ and $W_n$,
\begin{align}
x\,W_1(x;t,\alpha)&=\alpha t\,W_1(x-1;t-1,\alpha),\quad s,t\in\mathbb{N},\n
x_i\,W_n({\bm x};s,{\bm \beta})
&=\beta_is\,W_n({\bm x}-{\bm e}_i;s-1,{\bm \beta}),\quad i=1,\ldots,n,
\label{sumform2}
\end{align}
in which ${\bm e}_i$ is the unit vector in $i$-direction.
Similar sum formulas were used in \cite{os39} to solve many  generalisations 
of the original single variable problem.
It is interesting to note that the above formulas are easily generalised
\begin{align*}
(-x)_k\,W_1(x;t,\alpha)&=\alpha^k (-t)_k\,W_1(x-k;t-k,\alpha), \\
(-x_i)_k\,W_n({\bm x};s,{\bm \beta})
&=\beta_i^k(-s)_k\,W_n({\bm x}-k{\bm e}_i;s-k,{\bm \beta}),\quad k=1,\ldots.
\end{align*}
But I do not use them in this paper.

\bigskip
By plugging in the following form
\begin{equation}
P_{|m|=1}({\bm x})=1+\frac1N\sum_{i=1}^na_ix_i,
\label{onesol1}
\end{equation}
into the left eigenvalue equation
\begin{equation*}
\sum_{{\bm x}\in\mathcal{X}}K_n^{(i)}({\bm x},{\bm y}; N,{\bm\alpha}, {\bm \beta})
P_{|m|=1}({\bm x})=\lambda P_{|m|=1}({\bm y}),
\end{equation*}
and using the sum formulas \eqref{sum1form}, \eqref{sum2form}, one obtains
the degree one left eigenpolynomials.
\paragraph{For type (1) Markov chain}
The left eigenvector equation reads
\begin{align*}
1+\frac1N\sum_{i=1}^na_i\left\{\beta_iN-\beta_i\sum_{j=1}^n\alpha_jy_j+\alpha_iy_i\right\}
=\lambda\left(1+\frac1N\sum_{i=1}a_iy_i\right).
\end{align*}
These give rise to
\begin{align}
1+\sum_{j=1}^n\beta_ja_j&=\lambda,
\label{lam1eq}\\
-\left(\sum_{j=1}^n\beta_ja_j\right)\alpha_i+\alpha_ia_i&=\lambda a_i,\quad i=1,\ldots,n.
\label{char1eq0}
\end{align}
For each root  $\lambda_j$  of the characteristic equation
\begin{equation}
Det\left(\lambda\,I_n-F^{(1)}({\bm \alpha},{\bm \beta})\right)=0,\quad 
F^{(1)}({\bm \alpha},{\bm \beta})_{i\,j}\eqdef-\alpha_i\beta_j+\alpha_i\delta_{i\,j},\quad i,j=1,\ldots,n,
\label{char1eq}
\end{equation}
the coefficients $\{a_i\}$ are determined
\begin{equation*}
a_{i\,j}=-\frac{\alpha_i(\lambda_j-1)}{\lambda_j-\alpha_i},\quad i, j=1,\ldots,n.
\end{equation*}
 The consistency of the above formula and \eqref{lam1eq} implies the following sum rule
 \begin{equation}
\sum_{i=1}^n\frac{\alpha_i\beta_i}{\lambda_j-\alpha_i}=-1\   
\Longleftarrow \sum_{i=1}^n\beta_iu_{i\,j}=1-\lambda_j\quad j=1,\ldots,n,
\label{sum1rule}
\end{equation}
as $\lambda_j<1$ \eqref{leigeq}.
Let us tentatively identify the above $j$-th solution as $\bm{m}=\bm{e}_j$ left eigenpolynomial
\begin{equation}
P_{\bm{e}_j}(\bm{x};{\bm u})
=1-\frac1N\sum_{i=1}^n\frac{\alpha_i(\lambda_j-1)}{\lambda_j-\alpha_i}x_i,\quad j=1,\ldots,n.
\label{e_j1sol}
\end{equation}
\paragraph{For type (2) Markov chain}The left eigenvector equation reads
\begin{align*}
1+\frac1N\sum_{i=1}^na_i\left\{\beta_i(N-|y|)+\alpha_iy_i\right\}
=\lambda\left(1+\frac1N\sum_{i=1}a_iy_i\right).
\end{align*}
These give rise to
\begin{align}
1+\sum_{j=1}^n\beta_ja_j&=\lambda,
\label{lam2eq}\\
-\left(\sum_{j=1}^n\beta_ja_j\right)+\alpha_ia_i&=\lambda a_i,\quad i=1,\ldots,n.
\label{char2eq0}
\end{align}
For each root  $\lambda_j$  of the characteristic equation
\begin{equation}
Det\left(\lambda\,I_n-F^{(2)}({\bm \alpha},{\bm \beta})\right)=0,\quad 
F^{(2)}({\bm \alpha},{\bm \beta})_{i\,j}\eqdef-\beta_j+\alpha_i\delta_{i\,j},\quad i,j=1,\ldots,n,
\label{char2eq}
\end{equation}
the coefficients $\{a_i\}$ are determined
\begin{equation*}
a_{i\,j}=-\frac{\lambda_j-1}{\lambda_j-\alpha_i},\quad i, j=1,\ldots,n.
\end{equation*}
 The consistency of the above formula and \eqref{lam2eq} implies the following sum rule
 \begin{equation}
\sum_{i=1}^n\frac{\beta_i}{\lambda_j-\alpha_i}=-1\   
\Longleftarrow \sum_{i=1}^n\beta_iu_{i\,j}=1-\lambda_j\quad j=1,\ldots,n,
\label{sum2rule}
\end{equation}
as $\lambda_j<1$ \eqref{leigeq}.
Let us tentatively identify the above $j$-th solution as $\bm{m}=\bm{e}_j$ left eigenpolynomial
\begin{equation}
P_{\bm{e}_j}(\bm{x};{\bm u})
=1-\frac1N\sum_{i=1}^n\frac{\lambda_j-1}{\lambda_j-\alpha_i}x_i,\quad j=1,\ldots,n.
\label{e_j2sol}
\end{equation}
\begin{theo}
\label{uformula}{\bf System parameters}\
By the knowledge of the degree one left eigenvectors, 
\begin{align}
\sum_{{\bm x}\in\mathcal{X}}
K_n^{(i)}({\bm x},{\bm y}; N,{\bm\alpha}, {\bm \beta})
P_{\bm{e}_j}(\bm{x};{\bm u})=\lambda_jP_{{\bm e}_j}(\bm{y};{\bm u}),\quad j=1,\ldots,n,\quad i=1,2,
\label{deg1eq}
\end{align}
the system parameters $\{u_{i\,j}\}$ are completely identified,
\begin{align}
&P_{\bm{e}_j}(\bm{x};{\bm u})=1-\frac1N\sum_{i=1}^nu_{i\,j}x_i, \n
\Longrightarrow\quad  & 
u_{i\,j}=\frac{\alpha_i(\lambda_j-1)}{\lambda_j-\alpha_i}\qquad  \text{\rm for type\ (1)}  
\label{u1def} \\[-4pt]
&\hspace{70mm}i,j=1,\ldots,n, \n[-4pt]
& u_{i\,j}=
\frac{\lambda_j-1}{\lambda_j-\alpha_i} \qquad   \quad \, \text{\rm for type\ (2)}  
\label{u2def}
\end{align}
and each eigenvalue $\lambda_j$ satisfies the sum rule \eqref{sum1rule}, \eqref{sum2rule}.
As the degree one left eigenvectors  constitute the eigenvectors 
of the real symmetric matrix \eqref{Tneig1} they are orthogonal with 
each other and with $1=P_{\bm{0}}(\bm{x};{\bm u})$;
\begin{align}
\sum_{{\bm x}\in\mathcal{X}}W_n({\bm x};N,{\bm \eta})P_{\bm{e}_j}(\bm{x};{\bm u})=0,\
\sum_{{\bm x}\in\mathcal{X}}W_n({\bm x};N,{\bm \eta})
P_{\bm{e}_j}(\bm{x};{\bm u})P_{\bm{e}_k}(\bm{x};{\bm u})=0, \ j\neq k=1,\ldots,n.
\label{P1ort}
\end{align}
These are the so-called the orthogonality conditions {\rm \cite{Gri1, mizu, gr3}}.
It should be stressed that the system parameters $\{u_{i\,j}\}$ are irrational functions of the original
system parameters $\{\alpha_i\}$ and $\{\beta_i\}$.
 \end{theo}

From the characteristic equations \eqref{char1eq}, \eqref{char2eq}, 
one knows that the product and sum of the eigenvalues 
are rational functions of the original
system parameters $\{\alpha_i\}$ and $\{\beta_i\}$.
\begin{rema}
\label{eigprod}
\begin{align}
\text{\rm type \ (1)}\quad
\prod_{i=1}^n\lambda_i&=\prod_{i=1}^n\alpha_i\cdot(1-|\beta|),\qquad \quad
\sum_{i=1}^n\lambda_i=\sum_{i=1}^n\alpha_i(1-\beta_i),
\label{lam1sum}\\
\text{\rm type \ (2)}\qquad \prod_{i=1}^n\lambda_i
&=\prod_{i=1}^n\alpha_i-\sum_{j=1}^n\alpha_1\alpha_2\cdots\beta_j\cdots\alpha_n,
\qquad \sum_{i=1}^n\lambda_i=\sum_{i=1}^n(\alpha_i-\beta_i).
\label{lam2sum}
\end{align}
\end{rema}
\begin{rema}
\label{n=1} For $n=1$,
\begin{align}
\text{\rm type\ (1)}\qquad  \lambda&=\alpha(1-\beta),\quad 
u=\frac{1-\alpha+\alpha\beta}{\beta}=\eta^{-1},
\label{none1lam}\\
\text{\rm type\ (2)}\qquad \lambda&=\alpha-\beta,\qquad
\ u=\frac{1-\alpha+\beta}{\beta}=\eta^{-1}.
\label{none2lam}
\end{align}
\end{rema}
\begin{prop}{\bf $\mathfrak{S}_n$ Symmetry}
\label{Snsymm}
The Rahman polynomials  $\{P_{\bm{m}}(\bm{x})\}$\! \eqref{Pm} are
invariant under the symmetric group $\mathfrak{S}_n$ acting on ${\bm m}$, 
due to the arbitrariness of the ordering of $n$ roots $\{\lambda_j\}$ 
of the characteristic equations \eqref{char1eq}, \eqref{char2eq} in the parameters $u_{i\,j}$ 
\eqref{u1def}, \eqref{u2def}.
\end{prop}

\subsubsection{Exceptional cases}
\label{sec:excep}
For generic parameters, the eigenvalues are non-degenerate and the vanishing of a denominator
of $u_{i\,j}$, $\lambda_j-\alpha_i$ does not happen.
But as seen in the multivariate Krawtchouk and Meixner cases \cite{I23,mKrawt,mMeix},
certain exceptional parameter configurations do exist.
For example, when $\alpha_i=\alpha_j$,  the characteristic matrix \eqref{char1eq}, \eqref{char2eq} 
with $\lambda=\alpha_i$, {\em i.e.} $\alpha_i I_n-F^{(1)}({\bm \alpha},{\bm \beta})$ has equal
$i$-th and $j$-th rows  $-\alpha_i\beta_k$,  and $\alpha_i I_n-F^{(2)}({\bm \alpha},{\bm \beta})$ has equal
$i$-th and $j$-th rows $-\beta_k$, $k=1,\ldots,n$.
Thus the determinant vanishes and $\alpha_i$ becomes one of the eigenvalues, 
causing some $u_{i\,k}$ singular.
\begin{rema}
\label{rema:excep}
The parameters $\{\alpha_i\}$ should be distinct so that the general form of 
$P_{\bm m}({\bm x};{\bm u})$ \eqref{Pm} to hold.
\end{rema}
\begin{rema}
\label{rema:hoarah}
The degenerates case $\alpha_1=\alpha_2$ of the bivariate Rahman polynomials 
were discussed in some detail in {\rm \cite{HR}}.
\end{rema}

\subsection{Determination of the general eigenvalues}
\label{sec:eigval}

Now the whole of the Rahman polynomials are identified \eqref{Pm} and \eqref{u1def}, \eqref{u2def} and their 
orthogonality \eqref{pmort} is guaranteed by the above \eqref{P1ort} \cite{Gri1, mizu, mKrawt}.
The next task is to verify  the left eigenvalue equations
\begin{equation}
\sum_{{\bm x}\in\mathcal{X}}K_n^{(i)}({\bm x},{\bm y}; N,{\bm\alpha}, {\bm \beta})
P_{\bm{m}}(\bm{x};{\bm u})=\mathcal{E}({\bm m})P_{\bm{m}}(\bm{y};{\bm u}),
\quad \forall{\bm m}\in\mathcal{X},\quad i=1,2.
\label{KnPmeig}
\end{equation}
For this,  the knowledge of  the general eigenvalues $\mathcal{E}(\bm{m})$ is necessary. 
Since $P_{\bm m}({\bm 0};{\bm u})=1$, it is obtained by setting ${\bm y}={\bm 0}$,
\begin{align}
\mathcal{E}({\bm m})&=
\sum_{{\bm x}\in\mathcal{X}}K_n^{(i)}({\bm x},{\bm 0};{\bm \alpha},{\bm \beta})
P_{\bm m}({\bm x};{\bm u})
=\sum_{{\bm x}\in\mathcal{X}}\binom{N}{{\bm x}}P_{\bm m}({\bm x};{\bm u})
\prod_{i=1}^n\beta_i^{x_i}\cdot(1-|\beta|)^{N-|x|}\n
&=(1-|\beta|)^N\sum_{{\bm x}\in\mathcal{X}}\binom{N}{{\bm x}}P_{\bm m}({\bm x};{\bm u})
\prod_{i=1}^n\left(\frac{\beta_i}{1-|\beta|}\right)^{x_i}\n
&=(1-|\beta|)^N\sum_{{\bm x}\in\mathcal{X}}\binom{N}{{\bm x}}P_{\bm x}({\bm m};{\bm u}^T)
\prod_{i=1}^n\left(\frac{\beta_i}{1-|\beta|}\right)^{x_i}.
\label{em1}
\end{align}
Here, the duality 
\begin{equation*}
P_{\bm m}({\bm x};{\bm u})=P_{\bm x}({\bm m};{\bm u}^T)
\end{equation*}
of the $(n+1,2n+2)$ hypergeometric function \eqref{Pm} is used. 
That is, exchanging ${\bm x}\leftrightarrow {\bm m}$ is realised by $u_{i\,j}\leftrightarrow u_{j\,i}$.
The summation in r.h.s. of \eqref{em1}
has the same structure as that of the generating function $G({\bm x};{\bm u},{\bm t})$ in \eqref{gen}
except for the exchanges ${\bm x}\leftrightarrow {\bm m}$ and $u_{i\,j}\to u_{j\,i}$, 
${\bm t}\to{\bm \beta}/(1-|\beta|)$.
This means
\begin{align*}
\mathcal{E}({\bm m})&=(1-|\beta|)^NG\bigl({\bm m};{\bm u}^T,{\bm \beta}/(1-|\beta|)\bigr)\n
&=(1-|\beta|)^N\prod_{i=1}^n\left(1+\frac{|\beta|}{1-|\beta|}\right)^{N-|m|}
\left(1+\sum_{j=1}^nb_{i\,j}^T\beta_j/(1-|\beta|)\right)^{m_i}\\
&=\prod_{i=1}^n\left(1-|\beta|+\sum_{j=1}^nb_{i\,j}^T\beta_j\right)^{m_i}
=\prod_{i=1}^n\left(1-\sum_{j=1}^nu_{j\,i}\beta_j\right)^{m_i}=\prod_{i=1}^n\lambda_i^{m_i}.
\end{align*}
The explicit expression of $u_{i\,j}$ \eqref{u1def}, \eqref{u2def} and the sum rule 
\eqref{sum1rule}, \eqref{sum2rule}  
\begin{align*}
\text{for type (1)}\quad 
1-\sum_{j=1}^nu_{j\,i}\beta_j&=1-\sum_{j=1}^n\frac{\alpha_j\beta_j(\lambda_i-1)}{\lambda_i-\alpha_j}
=1+(\lambda_i-1)=\lambda_i,\\
\text{for type (2)}\quad 
1-\sum_{j=1}^nu_{j\,i}\beta_j&=1-\sum_{j=1}^n\frac{\beta_j(\lambda_i-1)}{\lambda_i-\alpha_j}
=1+(\lambda_i-1)=\lambda_i,
\end{align*}
lead to 
\begin{prop}
\label{prop:emformula} 
If $P_{\bm m}({\bm x};{\bm u})$ is the left eigenvector of 
$K_n^{(i)}({\bm x},{\bm y}; N,{\bm\alpha}, {\bm \beta})$ \eqref{KnPmeig}, it has a multiplicative spectrum
\begin{equation}
\mathcal{E}({\bm m})=\prod_{i=1}^n\lambda_i^{m_i}.
\label{emformula}
\end{equation}
\end{prop}

\subsection{Left eigenpvalue equation for $P_{\bm m}({\bm x};{\bm u})$}
\label{sec:lefteigeq}
The next and the main task is to prove the left eigenvalue equations for all the Rahman polynomials
\begin{equation}
\sum_{{\bm x}\in\mathcal{X}}K_n^{(i)}({\bm x},{\bm y}; N,{\bm\alpha}, {\bm \beta})
P_{\bm{m}}(\bm{x};{\bm u})=P_{\bm{m}}(\bm{y};{\bm u})\cdot\prod_{i=1}^n\lambda_i^{m_i},
\quad \forall{\bm m}\in\mathcal{X},\quad i=1,2.
\label{KnPmeig2}
\end{equation}
Let us multiply  both sides by $\binom{N}{\bm m}\prod_{i=1}^nt_i^{m_i}$ 
and take $\bm m$ summation to obtain
\begin{align*}
\sum_{{\bm x}\in\mathcal{X}}K_n^{(i)}({\bm x},{\bm y}; N,{\bm\alpha}, {\bm \beta})
\sum_{{\bm m}\in\mathcal{X}}\binom{N}{\bm m}P_{\bm{m}}(\bm{x};{\bm u})\prod_{k=1}^nt_k^{m_k}
&=\sum_{{\bm m}\in\mathcal{X}}\binom{N}{\bm m}P_{\bm{m}}(\bm{y};{\bm u})
\cdot\prod_{k=1}^n(\lambda_kt_k)^{m_k}.
\end{align*}
By comparing with the expansion formula \eqref{gen} of the generating function
\begin{equation}
G({\bm x};{\bm u},{\bm t})=(1+|t|)^{N-|x|}\prod_{i=1}^nT_i({\bm t})^{x_i},
\quad T_i({\bm t})\eqdef1+\sum_{j=1}^nb_{i\,j}t_j,
\end{equation}
the above equation reads
\begin{align}
\sum_{{\bm x}\in\mathcal{X}}
K_n^{(i)}({\bm x},{\bm y}; N,{\bm \alpha}, {\bm \beta})G({\bm x};{\bm u},{\bm t})
=G({\bm y};{\bm u},{\bm \lambda}{\bm t}),\quad i=1,2.
\label{desres}
\end{align}
The goal is achieved by proving this equation, which is rather tedious but straightforward.
\paragraph{For type (1) Markov chain}
\begin{align}
&\sum_{{\bm x}\in\mathcal{X}}
K_n^{(1)}({\bm x},{\bm y}; N,{\bm \alpha}, {\bm \beta})G({\bm x};{\bm u},{\bm t})\n
&=\sum_{{\bm x},{\bm z}\in\mathcal{X}}\binom{N-|z|}{{\bm x}-{\bm z}}\prod_{i=1}^n
\beta_i^{x_i-z_i}\left((1-|\beta|)(1+|t|)\right)^{N-|x|}T_i({\bm t})^{x_i}
\cdot\prod_{i=1}\binom{y_i}{z_i}\alpha_i^{z_i}(1-\alpha_i)^{y_i-z_i}\n
&=\sum_{{\bm x},{\bm z}\in\mathcal{X}}\binom{N-|z|}{{\bm x}-{\bm z}}\!\prod_{i=1}^n
\left(\beta_iT_i({\bm t})\right)^{x_i-z_i}\!\left((1-|\beta|)(1+|t|)\right)^{N-|x|}
\cdot
\prod_{i=1}\binom{y_i}{z_i}\!\left(\alpha_iT_i({\bm t})\right)^{z_i}(1-\alpha_i)^{y_i-z_i}\n
&=\sum_{{\bm z}\in\mathcal{X}}\Bigl((1-|\beta|)(1+|t|)
+\sum_{i=1}^n\beta_iT_i({\bm t})\Bigr)^{N-|z|}\cdot
\prod_{i=1}\binom{y_i}{z_i}\left(\alpha_iT_i({\bm t})\right)^{z_i}(1-\alpha_i)^{y_i-z_i}\n
&=(1+|\lambda t|)^N\prod_{i=1}^n
\sum_{z_i=1}^{y_i}\binom{y_i}{z_i}
\left(\frac{\alpha_iT_i({\bm t})}{1+|\lambda t|}\right)^{z_i}(1-\alpha_i)^{y_i-z_i}.
\label{1srop}
\end{align}
Here the simplification of the first term is used.
\begin{align}
\sum_{i=1}^n\beta_iT_i({\bm t})&=\sum_{i=1}^n\beta_i\left(1+\sum_{j=1}^n(1-u_{i\,j})t_j\right)
=|\beta|(1+|t|)-\sum_{i\,j}\beta_iu_{i\,j}t_j\n
&=|\beta|(1+|t|)-|t|+|\lambda t|
 \Longrightarrow (1-|\beta|)(1+|t|)+\sum_{i=1}^n\beta_iT_i({\bm t})=1+|\lambda t|.
\label{simp2}
\end{align}
At the last equality, the sum rule \eqref{sum1rule} is used. 
Now, the $z_i$ summation gives
\begin{align*}
\sum_{z_i=1}^{y_i}\binom{y_i}{z_i}
\left(\frac{\alpha_iT_i({\bm t})}{1+|\lambda t|}\right)^{z_i}(1-\alpha_i)^{y_i-z_i}
&=\left(1-\alpha_i+\frac{\alpha_iT_i({\bm t})}{1+|\lambda t|}\right)^{y_i}\\
&=(1+|\lambda t|)^{-y_i}\bigl ((1-\alpha_i)(1+|\lambda t|)+\alpha_iT_i({\bm t})\bigr).
\end{align*}
A straightforward calculation based on the explicit form of $u_{i\,j}$ \eqref{u1def} shows
\begin{equation}
(1-\alpha_i)(1+|\lambda t|)+\alpha_iT_i({\bm t})=T_i({\bm \lambda}{\bm t}).
\label{straight}
\end{equation}
This leads to the desired result \eqref{desres}.
\paragraph{For type (2) Markov chain}
\begin{align}
&\sum_{{\bm x}\in\mathcal{X}}
K_n^{(2)}({\bm x},{\bm y}; N,{\bm \alpha}, {\bm \beta})G({\bm x};{\bm u},{\bm t})\n
&=\sum_{{\bm x},{\bm z}\in\mathcal{X}}\binom{N-|y|}{{\bm x}-{\bm z}}\prod_{i=1}^n
\beta_i^{x_i-z_i}(1-|\beta|)^{N-|x|-|y|+|z|}\cdot (1+|t|)^{N-|x|}T_i({\bm t})^{x_i}\n
&\hspace{50mm}
\times\prod_{i=1}\binom{y_i}{z_i}\alpha_i^{z_i}(1-\alpha_i)^{y_i-z_i}\n
&=\sum_{{\bm x},{\bm z}\in\mathcal{X}}\binom{N-|y|}{{\bm x}-{\bm z}}\!\prod_{i=1}^n
\left(\beta_iT_i({\bm t})\right)^{x_i-z_i}\!\left((1-|\beta|)(1+|t|)\right)^{N-|x|-|y|+|z|}\n
&\hspace{50mm}\times
\prod_{i=1}\binom{y_i}{z_i}\!\left(\alpha_iT_i({\bm t})\right)^{z_i}\bigl((1+|t|)(1-\alpha_i)\bigr)^{y_i-z_i}\n
&=\Bigl((1-|\beta|)(1+|t|)
+\sum_{i=1}^n\beta_iT_i({\bm t})\Bigr)^{N-|y|}\n
&\hspace{50mm}
\times\sum_{{\bm z}\in\mathcal{X}}\prod_{i=1}\binom{y_i}{z_i}\left(\alpha_iT_i({\bm t})\right)^{z_i}\bigl((1+|t|)(1-\alpha_i)\bigr)^{y_i-z_i}\n
&=(1+|\lambda t|)^{N-|y|}\prod_{i=1}^n
\Bigl((1+|t|)(1-\alpha_i)+\alpha_iT_i({\bm t})\Bigr)^{y_i}.
\label{2srop}
\end{align}
In the first term \eqref{simp2} is used and  for the second term 
\eqref{straight} is used, which also holds for type (2).
The desired result \eqref{desres} is obtained.

\begin{theo}
\label{theo:main}
The Rahman polynomials $\{P_{\bm m}({\bm x};{\bm u})\}$ \eqref{Pm} are the left eigenvectors 
of the Markov chain operator $K_n^{(i)}({\bm x},{\bm y}; N,{\bm\alpha}, {\bm \beta})$ with the
multiplicative eigenvalues $\prod_{i=1}^n\lambda_i^{m_i}$ \eqref{emformula}, in which $\lambda_i$'s are
the eigenvalues of the degree one polynomials,   i.e. 
the roots of the $n\times n$ characteristic equation \eqref{char1eq}, \eqref{char2eq},
\begin{equation*}
\sum_{{\bm x}\in\mathcal{X}}K_n^{(i)}({\bm x},{\bm y}; N,{\bm\alpha}, {\bm \beta})
P_{\bm{m}}(\bm{x};{\bm u})=P_{\bm{m}}(\bm{y};{\bm u})\cdot\prod_{i=1}^n\lambda_i^{m_i},
\quad \forall{\bm m}\in\mathcal{X},\quad i=1,2.
\tag{\ref{KnPmeig2}}
\end{equation*}
Likewise $\{P_{\bm m}({\bm x};{\bm u})W_n({\bm x};N,{\bm \eta})\}$  are the eigenvectors of 
$K_n^{(i)}({\bm x},{\bm y}; N,{\bm\alpha}, {\bm \beta})$ with the same
 eigenvalue $\prod_{i=1}^n\lambda_i^{m_i}$ \eqref{emformula}, 
 \begin{align}
\sum_{{\bm y}\in\mathcal{X}}K_n^{(i)}({\bm x},{\bm y}; N,{\bm\alpha}, {\bm \beta})
P_{\bm{m}}(\bm{y};{\bm u})W_n({\bm y};N,{\bm \eta})=\prod_{i=1}^n\lambda_i^{m_i}\cdot 
P_{\bm{m}}(\bm{x};{\bm u})W_n({\bm x};N,{\bm \eta}),\n
\quad \forall{\bm m}\in\mathcal{X},\quad i=1,2,
\label{KnPmeig3}
\end{align}
in which $W_n({\bm x};N,{\bm \eta})$ is the reversible distribution for 
$K_n^{(i)}({\bm x},{\bm y}; N,{\bm\alpha}, {\bm \beta})$ \eqref{Kn1rev}, \eqref{Kn2rev}.
\end{theo}

%
%
\section{Comments}
\label{sec:comm}

The most important point of this paper is that the Rahman polynomials are obtained as the
{\bf left eigenvectors }of the Markov chain operators 
$K_n^{(i)}({\bm x},{\bm y}; N,{\bm\alpha}, {\bm \beta})$.
Therefore, they are completely determined by the initial system parameters $\{{\bm \alpha}\}$
and $\{{\bm \beta}\}$ alone  as shown in {\bf Theorems \ref{uformula} \ref{theo:main}}.
In contrast, in some literature, the {\bf eigenvalue problem} was addressed \cite{HR0}(3.17), 
\cite{HR}(1.5), \cite{gr1}\S6, \cite{gr3}(1.7). This naturally involves the reversible distribution
$W_n({\bm x};N,{\bm \eta})$
 \begin{equation*}
\sum_{{\bm y}\in\mathcal{X}}K_n^{(i)}({\bm x},{\bm y}; N,{\bm\alpha}, {\bm \beta})
P_{\bm{m}}(\bm{y};{\bm u})W_n({\bm y};N,{\bm \eta})=\prod_{i=1}^n\lambda_i^{m_i}\cdot 
P_{\bm{m}}(\bm{x};{\bm u})W_n({\bm x};N,{\bm \eta}),
\quad \forall{\bm m}\in\mathcal{X},
\end{equation*}
and the parameters ${\bm \eta}$ come into the play. 
Efforts were made to construct Rahman polynomials from the fact that they had the orthogonality 
measure $W_n({\bm x};N,{\bm \eta})$ and ${\bm \eta}$'s came into the polynomial parameters.
Of course some combinations of the original parameters ${\bm \alpha}$ and ${\bm \beta}$
can be expressed by ${\bm \eta}$'s based on the relation I.(1.6).
Therefore, it might be possible to mingle some ${\bm \eta}$'s 
into the expression of the system parameters $u_{i\,j}$. 
But it is totally unnecessary.

To say symbolically,
``The parameters ${\bm \eta}$'s are chosen by 
$K_n^{(i)}({\bm x},{\bm y}; N,{\bm\alpha}, {\bm \beta})$ to provide the orthogonality weight 
for the Rahman polynomials.
The orthogonality weights never determine multivariate orthogonal polynomials,
as shown explicitly for the cases of the Krawtchouk and Meixner polynomials \cite{mKrawt,mMeix}."

In some works, the importance of the left eigenvalue problems was well recognised
\cite{coo-hoa-rah77}(2.2), \cite{HR0}(3.23), \cite{gr3}\S6, but unfortunately 
the idea was not well developed.
A superficial reason might be that they adopted notation like 
$K_n({\bm j},{\bm i}; N,{\bm\alpha}, {\bm \beta})$ \cite{gr3}(1.7)
which was rather distant from polynomials in ${\bm x}$, 
$P_{\bm m}({\bm x};{\bm u})$ as compared with mine, 
$K_n^{(i)}({\bm x},{\bm y}; N,{\bm\alpha}, {\bm \beta})$.
A similar situation existed in the treatments of birth and death processes.
They used $\lambda_n$ and $\mu_n$ for the birth/death rates \cite{KarMcG} instead of $B(x)$ and $D(x)$
\cite{bdsol} and $B_j({\bm x})$ and $D_j({\bm x})$ for the multivariate cases \cite{mKrawt,mMeix}.

It should be emphasised that the eigenvalues were correctly identified \cite{HR, gr3} in terms of $u_{i\,j}$,
as the eigenvalues are common in the eigenvalue problem and the left eigenvalue problem.
For the bivariate case \cite{HR},
\begin{equation*}
\lambda_{m,n}=(1-\beta_1t-\beta_2u)^m(1-\beta_1v-\beta_2w)^n,
\tag{\cite{HR}.4.23}
\end{equation*}
Which is the same as the sum rule \eqref{sum1rule} and $\lambda_{m,n}=\lambda_1^m\lambda_2^n$,
as 
$
\left(
\begin{array}{cc}
t  &   v   \\
u  &  w       
\end{array}
\right)
=
\left(
\begin{array}{cc}
u_{1\,1}  &   u_{1\,2}   \\
u_{2\,1} &  u_{2\,2}     
\end{array}
\right)
$.
For the general multivariate case \cite{gr3}\S7, $\lambda_{\bm m}=\prod_{i=1}^n(1-\omega_i)^{m_i}$,
in which $\omega_i=\sum_{j=1}\beta_ju_{i\,j}$. This is again the same as 
my multiplicative spectrum \eqref{emformula} by using the sum rule \eqref{sum1rule}. 
Their $u_{i\,j}$ is my $u_{j\,i}$.
The explicit expressions of $u_{i\,j}$ in terms of the original system parameters 
${\bm \alpha}$ and ${\bm \beta}$ 
were not available in their works.

In \cite{os39}, various single variable discrete time Markov chains were constructed based on
five types of convolutions of the orthogonality measures $\pi(x, N,{\bm \lambda})$ of
the Krawtchouk, Charlier, Hahn, Meixner, $q$-Hahn and  $q$-Meixner polynomials 
\cite{askey, ismail, koeswart}. The five types of the convolutions were
\begin{align}
  \text{(\romannumeral1)}:&\ \ K(x,y)\eqdef\sum_{z=0}^{\min(x,y)}
  \!\pi(x-z,N-z,\bm{\lambda}_2)\pi(z,y,\bm{\lambda}_1),
  \label{conv11}\\
  \text{(\romannumeral2)}:&\ \ K(x,y)\eqdef\sum_{z=\max(0,x+y-N)}^{\min(x,y)}
  \!\!\!\!\!\!\!\!\!\pi(x-z,N-y,\bm{\lambda}_2)\pi(z,y,\bm{\lambda}_1),
  \label{conv2}\\
  \text{(\romannumeral3)}:&\ \ K(x,y)\eqdef\sum_{z=\max(x,y)}^N
  \!\!\!\!\pi(x,z,\bm{\lambda}_2)\pi(z-y,N-y,\bm{\lambda}_1),
  \label{conv3}\\
  \text{(\romannumeral4)}:&\ \ K(x,y)\eqdef\sum_{z_2=0}^{\min(x,y)}
  \!\!\pi(z_2,y,\bm{\lambda}_1)\!\!\!\sum_{z_1=\max(x,y)}^N
  \!\!\!\!\pi(x-z_2,z_1-z_2,\bm{\lambda}_3)\pi(z_1-y,N-y,\bm{\lambda}_2),
  \label{conv4}\\
  \text{(\romannumeral5)}:&\ \ K(x,y)\eqdef\sum_{z_2=0}^{\min(x,y)}
  \!\!\pi(z_2,y,\bm{\lambda}_1)\!\!\!\sum_{z_1=x+y-z_2}^N
  \!\!\!\!\pi(x-z_2,z_1-y,\bm{\lambda}_3)\pi(z_1-y,N-y,\bm{\lambda}_2).
  \label{conv5}
\end{align}
In these formulas ${\bm \lambda}_i$, $i=1,2,3$ was a set of parameters.
 These $K(x,y)$  all provided exactly solvable single variable Markov chains.
The multivariate versions of type (i)  and (ii) convolutions  correspond to $K_n^{(1)}$
and   $K_n^{(2)}$, respectively.
I wonder if another Rahman like polynomials could be obtained by 
certain multivariable generalisations of these convolutions.


\goodbreak
\end{document}